\documentclass[12pt,a4paper,reqno]{amsart}

\makeatletter
\def\author@andify{%
  \nxandlist {\unskip ,\penalty-1 \space\ignorespaces}%
    {\unskip {} \@@and~}%
    {\unskip \penalty-2 \space \@@and~}%
}
\makeatother

\usepackage[english]{babel}
\usepackage{amssymb,latexsym,amsfonts,amsthm,upref,amsmath}
\usepackage[margin=1in]{geometry} 
\usepackage[dvipsnames,x11names]{xcolor}
 \definecolor{myblue}{HTML}{003399}
\usepackage{hyperref}
\hypersetup{colorlinks,citecolor=myblue,filecolor=black,linkcolor=myblue,urlcolor=myblue}
\usepackage{enumerate}  
\usepackage{comment}  
\usepackage{tikz}
\usepackage{float}
\usepackage{cleveref}
\usepackage{mathtools}
\usepackage{cite}
\usepackage{filecontents}
\usepackage{enumitem}
\makeatletter
\newcommand{\leqnomode}{\tagsleft@true}
\newcommand{\reqnomode}{\tagsleft@false}

\makeatother
\newtheorem*{thm*}{Theorem}
\newtheorem*{lem*}{Lemma}
\newtheoremstyle{prim}{}{}{\normalfont}{}{\bfseries}{.}{ }{}
\newtheoremstyle{stil}{}{}{\slshape}{}{\bfseries}{.}{ }{}
\theoremstyle{stil}
\newtheorem{thm}{Theorem}[section]
\newtheoremstyle{defi}{}{}{}{}{\bfseries}{.}{ }{}
\theoremstyle{defi}

\theoremstyle{defi}
\newtheorem{rem}[thm]{Remark}
\theoremstyle{stil}
\newtheorem*{mthm*}{Main Theorem}
\newtheorem*{kor*}{Corollary}
\newtheorem{pro}[thm]{Proposition}
\theoremstyle{stil}
\newtheorem{lem}[thm]{Lemma}
\theoremstyle{stil}
\newtheorem{kor}[thm]{Corollary}
\theoremstyle{prim}

\newenvironment{prf}{\noindent \textit{Proof.}}{\null\hfill$\qed$\hskip
2mm\vskip 2mm}

\newcommand{\dygl}{ \text{DY}(\mathfrak{gl}_2)}
\newcommand{\dysl}{ \text{DY}(\mathfrak{sl}_2)}
\newcommand{\gl}{\mathfrak{gl}_2}
\newcommand{\sll}{\mathfrak{sl}_2}
\newcommand{\h}{\mathfrak{h}}
\newcommand{\s}{\mathfrak{s}}

\newcommand{\vac}{\mathop{\mathrm{\boldsymbol{1}}}}

\newcommand{\CC}{\mathbb{C}}

\newcommand{\ZZ}{\mathbb{Z}}

\newcommand{\norm}[1]{\hspace{3pt}:\mathrel{#1}:\hspace{3pt}}

\newcommand{\dyg}{ \text{DY}(\mathfrak{gl}_2) }
\newcommand{\dys}{ \text{DY}(\mathfrak{sl}_2) }
\newcommand{\hh}{\text{H}}

\newcommand{\wtld}{\widetilde}
\newcommand{\wht}{\widehat}
\newcommand{\wvr}{\overline}

\newcommand{\ot}{\otimes}
\newcommand{\ts}{\hspace{1pt}}
\newcommand{\qdet}{ {\rm qdet}\hspace{1pt}}

\newcommand{\ndo}{\mathop{\mathrm{End}}}
\newcommand{\om}{\mathop{\mathrm{Hom}}}

\newcommand{\fand}{\quad\text{and}\quad}
\newcommand{\Fand}{\qquad\text{and}\qquad}

\newcommand{\non}{\nonumber}
\newcommand{\beq}{\begin{equation}}
\newcommand{\eeq}{\end{equation}}
\newcommand{\ben}{\begin{equation*}}
\newcommand{\een}{\end{equation*}}

\makeatletter
\def\smalloverbrace#1{\mathop{\vbox{\m@th\ialign{##\crcr\noalign{\kern3\p@}%
  \tiny\downbracefill\crcr\noalign{\kern3\p@\nointerlineskip}%
  $\hfil\displaystyle{#1}\hfil$\crcr}}}\limits}
\makeatother

\makeatletter
\def\smallunderbrace#1{\mathop{\vtop{\m@th\ialign{##\crcr
   $\hfil\displaystyle{#1}\hfil$\crcr
   \noalign{\kern3\p@\nointerlineskip}%
   \tiny\upbracefill\crcr\noalign{\kern3\p@}}}}\limits}
\makeatother

\begin{document}

\title{Semi-infinite construction for the double Yangian of type $A_1^{(1)}$}

\author{Marijana Butorac}
\address{Marijana Butorac:\newline\indent
Faculty of Mathematics, University of Rijeka,\newline\indent Radmile Matej\v{c}i\'{c} 2, 51000 Rijeka, Croatia}
\email{mbutorac@math.uniri.hr}

\author{Naihuan Jing}
\address{Naihuan Jing:\newline\indent Department of Mathematics, North Carolina State University,\newline\indent  Raleigh, NC 27695, USA}
\email{jing@ncsu.edu}

\author{Slaven Ko\v{z}i\'{c}} 
\address{Slaven Ko\v{z}i\'{c}:\newline\indent  Department of Mathematics, Faculty of Science, University of Zagreb,\newline\indent   Bijeni\v{c}ka cesta 30, 10000 Zagreb, Croatia}
\email{kslaven@math.hr}

\author{Fan Yang} 
\address{Fan Yang:\newline\indent  Department of Mathematics, Jiaying University,\newline\indent    Meizhou, Guangdong 514000, China}
\email{1329491781@qq.com}


\begin{abstract}
We  consider   certain   infinite dimensional modules of level 1 for  the double Yangian $\text{DY}(\mathfrak{gl}_2)$
which are based on the   Iohara--Kohno   realization. We show that they possess    topological bases of Feigin--Stoyanovsky-type, i.e. the bases expressed in terms   of semi-infinite   monomials of    certain   integrable operators    which stabilize and satisfy the difference two condition.
Finally, we give some applications of these bases to the representation theory of the corresponding quantum affine vertex algebra.
\end{abstract}

\maketitle

\allowdisplaybreaks


\section{Introduction}\label{intro}
\numberwithin{equation}{section}

 The integrable
highest weight modules present one of the most fundamental notions in the  representation theory of affine Kac--Moody Lie algebras;  
 see, e.g., the book by Kac \cite{K}. The problem of constructing   different types of  bases for such modules and their various substructures, especially those which establish connection with  Rogers--Ramanujan-type identities via  character formulae, has been extensively studied since the pioneering paper of Lepowsky and Milne \cite{LM}.
Our paper is motivated by  the  well known   Feigin--Stoyanovsky  construction \cite{FS} of semi-infinite monomial bases for  certain integrable highest weight modules for the affine Lie algebra $\widehat{\mathfrak{sl}}_2$. The  construction relies on the fact that these modules
 can be obtained from their   principal subspaces  using the Weyl translation operator. At the level 1, the resulting   bases consist  of   semi-infinite monomials  $x_\alpha(r_1)x_\alpha(r_2)\ldots$ in coefficients of the vertex operator $x_\alpha (z)  =\sum_{r\in\ZZ} x_\alpha(r)z^{-r-1}$ associated with the positive simple root $\alpha$ of $\mathfrak{sl}_2$. Their degrees $r_1,r_2,\ldots $ satisfy the   difference two condition  $ r_{j+1}\geqslant r_j+2$ for all $j=1,2,\ldots ,$  which comes from the   integrability  relation  $x_\alpha (z)^2=0$ of Lepowsky and Primc \cite{LP}. Moreover, these monomials   stabilize, i.e. for a sufficiently large  $n$ all degrees $r_n, r_{n+1},\ldots$ are consecutive odd or even integers, depending on the choice of the highest weight. Later on, the semi-infinite construction was generalized to the case of quantum affine algebra $U_q(\widehat{\mathfrak{sl}}_2)$ by Ding and B. Feigin \cite{DF} using the realization of its  integrable highest weight modules found by I. Frenkel and the second author in \cite{FJ}. 

The goal of this paper is to give a semi-infinite construction for certain infinite dimensional modules  of level 1  for the centrally extended double Yangian $\text{DY}(\mathfrak{gl}_2)$    defined over the commutative ring $\CC[[h]]$. Their bosonic realization, which resembles the famous Frenkel--Kac--Segal construction \cite{FK,S} for affine Lie algebras, was   given by Iohara and Kohno in \cite{IK} for $\text{DY}(\mathfrak{gl}_2)$ and then generalized to the higher rank case by Iohara \cite{I}.  We slightly  modify   the   Iohara--Kohno realization as the  action of the original translation operator   does not appear to be in tune with the semi-infinite construction. However, the action of the double Yangian   is still given on the same $\CC[[h]]$-module, which we denote by $\mathcal{F}_i$, $i=0,1$. In contrast with the aforementioned setting of affine Lie algebras and quantum affine algebras,   the general theory 
of integrable representations for double Yangians has not  yet  been sufficiently developed. Thus, in our construction,   we often need to use different and  more technical arguments which rely  on the explicit formulae for the action of the double Yangian generators on $\mathcal{F}_i$. 

Motivated by  Ding--Feigin's approach \cite{DF}, we start by defining an auxiliary commutative operator   $\wvr{X}(z)$ on  $\mathcal{F}_i$, $i=0,1$, which can be regarded as a  Yangian counterpart  of  the level 1 affine vertex operator $x_\alpha (z) $. In particular, it satisfies the   $h$-adic integrability  relation  $\wvr{X}(z)\wvr{X}(z\pm h)=0$. We use its coefficients  in parallel with \cite{FS,G} to introduce   the notion of   principal submodule  $W_i \subset \mathcal{F}_i$ and, furthermore,  to obtain the   topological 
basis for $W_i$ which provides an interpretation of the sum-sides of Rogers--Ramanujan identities.  
Next, we employ  the action of translation operator on $W_i$  to recover     irreducible  modules $\mathcal{L}_i(\mathfrak{sl}_2)\subset \mathcal{F}_i$ for the double Yangian $\text{DY}(\mathfrak{sl}_2)$,  such that their classical limits   are exactly the level $1$ integrable  highest weight  $\widehat{\mathfrak{sl}}_2$-modules $L(\Lambda_i)$ of highest weight $\Lambda_i$ with $i=0,1$. Finally, we construct  the Feigin--Stoyanovsky-type semi-infinite monomial bases    for $\mathcal{L}_i(\mathfrak{sl}_2)$, which is the main result of this paper. In addition, we generalize this construction   to the corresponding  modules $\mathcal{L}_i(\mathfrak{gl}_2)$ for the double Yangian $\text{DY}(\mathfrak{gl}_2)$ by using the action of its Heisenberg subalgebra,  which is generated by the coefficients of the quantum determinant and commutes with $\text{DY}(\mathfrak{sl}_2)$.

At the end of the paper, we  obtain some applications of the semi-infinite construction to the quantum vertex algebra theory. In particular, by employing the Iohara--Kohno isomorphism \cite{IK} between two realizations of the double Yangian, we show that $\mathcal{L}_i(\mathfrak{gl}_2)$ are naturally equipped with the structure of   irreducible modules for the corresponding Etingof--Kazhdan  quantum affine vertex algebra of level 1 from \cite{EK}.

\section{Preliminaries}\label{sec0}

In this section, we 
recall the double Yangian for the general linear  Lie algebra  $\mathfrak{gl}_2$ and the  Iohara--Kohno bosonic realization of its level 1 modules.

\subsection{Double Yangian for $\mathfrak{gl}_2$ }\label{subsec00}
\setcounter{equation}{0}

We follow the paper of Iohara and Kohno \cite{IK} to introduce the centrally extended double Yangians for the Lie algebras  $\mathfrak{gl}_2$ and $\mathfrak{sl}_2$ and recover some of their properties.  Let $I$ be the identity and $P$ the permutation operator on $ \CC^2\ot \CC^2$. Consider the (normalized) {\em Yang $R$-matrix} over the commutative ring $\CC[[h]]$, 
$$
R(u)=\frac{1}{1+h/u}\left(I +\textstyle \frac{h}{u}P\right)\in \ndo\CC^2\ot\ndo\CC^2[[h/u]].
$$

The {\em double Yangian} $\dyg$  is defined as the    associative algebra over the   ring $\CC[[h]]$ generated by the central element $C$ and the elements $t_{ij}^{(  r)}$, where $i,j=1,2$ and $r\in\ZZ.$ Its defining relations are given by
\begin{align}
R_{12}(u-v)\ts T_{13}^{\pm}(u)\ts T_{23}^{\pm}(v) &
=T_{23}^{\pm}(v)\ts T_{13}^{\pm}(u)\ts R_{12}(u-v),\label{rtt1}\\
R_{12}(u-v-hC/2)\ts T_{13}^{+}(u)\ts T_{23}^{-}(v)  &
=T_{23}^{-}(v)\ts T_{13}^{+}(u)\ts R_{12}(u-v+hC/2).\label{rtt2}
\end{align}
The generator matrices $T^{\pm}(u)$ are defined by
\beq\label{tu}
T^{\pm}(u)=\sum_{i,j=1,2} e_{ij}\ot t_{ij}^{\pm} (u),
\eeq
where   $e_{ij}\in \ndo\CC^2$ denote the matrix units and the power series $t_{ij}^{\pm}(u)$ are given by
\beq\label{tij}
t_{ij}^+(u)=\delta_{ij}-h\sum_{r\geqslant 0}t_{ij}^{(r)}\ts u^{-r-1}
\Fand
t_{ij}^-(u)=\delta_{ij}+h\sum_{r\geqslant 1}t_{ij}^{(-r)}\ts u^{r-1}
.
\eeq
In    \eqref{rtt1} and \eqref{rtt2} we use the subscripts to indicate the   tensor factors, i.e. we have
$$
T_{13}^{\pm}(u)=\sum_{i,j=1,2} e_{ij}\ot 1\ot t_{ij}^{\pm} (u)
\Fand
T_{23}^{\pm}(u)=\sum_{i,j=1,2} 1\ot e_{ij} \ot t_{ij}^{\pm} (u).
$$

Let us discuss the classical limit of   double Yangian. Consider  the affine Lie algebra 
$\wht{\mathfrak{gl}}_2=\mathfrak{gl}_2\ot\CC[ t^{\pm 1}]\oplus \CC K $,  where $K$ is the   central element and the Lie brackets are  
\beq\label{commrel}
\left[e_{ij}(r),e_{kl}(s)\right]
=\delta_{kj}\ts e_{il} (r+s) -\delta_{il}\ts e_{kj}(r+s)
+r\ts \delta_{r+s\ts 0}\ts K\left(\delta_{kj}\ts\delta_{il} -\delta_{ij}\ts \delta_{kl}\right) 
\eeq
for 
$e_{ij}(r)=e_{ij}\ot t^r$. Introduce the ascending filtration over 
$\dyg$ by setting
\beq\label{degree}
\deg t_{ij}^{(r)}=r
\quad\text{for}\quad i,j=1,2,\,r\in\ZZ
\Fand \deg C= 0.
\eeq 
The images  $\bar{t}_{ij}^{(r)}$ and $\bar{C}$ of the double Yangian generators $ t_{ij}^{(r)}$ and $C$ in the corresponding   graded algebra $\text{gr}\ts \dyg$ satisfy    \eqref{commrel}. Thus, the assignments
$  
e_{ij}(r) \mapsto \bar{t}_{ij}^{(r)}$ and  $K \mapsto \bar{C}
$
define the algebra homomorphism
\beq\label{hom}
U(\wht{\mathfrak{gl}}_2)\ot\CC[[h]]  \to \text{gr}\ts \dyg
.
\eeq

\begin{thm}
The map \eqref{hom} is an algebra isomorphism.
\end{thm}

\begin{prf}
The surjectivity of the map \eqref{hom} is clear. On the other hand, the injectivity is a consequence of the Poincar\'{e}--Birkhoff--Witt theorem \cite[Thm. 2.2]{JKMY}, which states that the suitably ordered monomials in the double Yangian generators form its basis; see also \cite[Prop. 3.1]{JY} and \cite[Thm. 15.3]{N}. More specifically, although we define the double Yangian using  
 the  normalization of  the Yang $R$-matrix which differs from   \cite{JKMY}, the arguments from the corresponding part of the proof of \cite[Thm. 2.2]{JKMY} can be still carried out analogously. It is worth noting that they rely on  the Iohara--Kohno realization   \cite{IK}, which provides level $1$ representations of the double Yangian; see Theorem \ref{IKthm} below.
\end{prf}

\begin{rem}\label{klasicnilimes}
The correspondence similar to \eqref{hom}, which employs the  universal  enveloping algebra over $\CC$,  can be also   established by taking the classical limit  $\dyg / h\dyg$ of the double Yangian. Indeed,  by extracting the coefficients of the matrix entries in its defining relations  \eqref{rtt1} and \eqref{rtt2},  one  observes that the resulting  top degree terms   with respect to \eqref{degree} coincide with the terms which contain the lowest power of   $h$. 
\end{rem}

From now on, we shall assume that the double Yangian for $\mathfrak{gl}_2$  is $h$-adically completed.
Using the series \eqref{tij} one    obtains  its Drinfeld generators    \cite{D} as follows:
\begin{align}
&k_1^{\pm}(u)=t_{11}^{\pm}(u),\qquad
k_2^{\pm}(u)=t_{22}^{\pm}(u) - t_{21}^{\pm}(u)\ts t_{11}^{\pm}(u)^{-1}\ts t_{12}^{\pm} (u),\label{rels1}\\
&X^+(u)=t_{11}^+(u-hC/4)^{-1}\ts t_{12}^+(u-hC/4)-t_{11}^-(u+hC/4)^{-1}\ts t_{12}^-(u+hC/4),\label{rels2}\\
&X^-(u)=t_{21}^+(u+hC/4)\ts t_{11}^+(u+hC/4)^{-1} - t_{21}^-(u-hC/4)\ts t_{11}^- (u-hC/4)^{-1}. \label{rels3}
\end{align}
The commutation relations for these generators can be found in \cite[Thm. 2.1]{IK}.

The double Yangian $\dyg$ can be decomposed into two subalgebras: the {\em double Yangian} $\dys$, which is generated by the   central element $C_1\coloneqq C$ and all coefficients of the power series
\begin{gather}
E(u)=\frac{1}{h}X^+(u+h/2),\quad F(u)=\frac{1}{h}X^-(u+h/2),\non\\ 
H^{\pm}(u)=k_2^{\pm }(u+h/2)\ts k_1^{\pm}(u+h/2)^{-1},\label{hplus}
\end{gather} 
and the {\em Heisenberg subalgebra} $\hh$, which is generated by the   coefficients of the  series
$$
K^{\pm}(u)=k_1^{\pm }(u-h/2)\ts k_2^{\pm}(u+h/2) 
$$
and the central element $C_2= -2C$.
The  generators of the Heisenberg subalgebra commute with all elements of $\dys$
and satisfy the   relations from \cite[Cor. 2.2]{IK},
\begin{gather*}
\left[K^{\pm}(u),K^{\pm}(v)\right]=0,\\
 \frac{u-v-h+hC_2/4}{u-v+h+hC_2/4}K^+(u)\ts K^-(v)=
K^-(v)\ts K^+(u) \frac{u-v-h-hC_2/4}{u-v+h-hC_2/4} .
\end{gather*}

The generator series of $\hh$   are of the form
$ 
K^{\pm}(u)=1\mp h\kappa^{\pm} (u)  
$,
where
$$
\kappa^+ (u)=\sum_{r\geqslant 0}\kappa^{(r)}u^{-r-1}\Fand 
\kappa^- (u)=\sum_{r\geqslant 1}\kappa^{(-r)}u^{r-1}.
$$
Consider the Heisenberg Lie subalgebra $\wht{\mathfrak{h}}=\mathfrak{h}\ot\CC[ t^{\pm1}]\ot\CC K_2\subset \wht{\mathfrak{gl}}_2$, where $K_2=-2K$ is the central element and $\mathfrak{h}=\CC I$ is the one-dimensional commutative subalgebra of $ \mathfrak{gl}_2$ spanned by the identity matrix $I=e_{11}+e_{22}$. By \eqref{commrel}, its Lie brackets are given by
$$
\left[I(r),I(s)\right]
= r\ts \delta_{r+s\ts 0}\ts K_2,
$$
where $I(r)=I\ot t^r$.
Note that the restriction of  the map \eqref{hom}   produces the   isomorphism 
\beq\label{hizo}
U(\wht{\mathfrak{h}})\ot\CC[[h]]\to \text{gr}\ts \hh,
\eeq
 where $\text{gr}\ts \hh$ is the  corresponding graded algebra of  $  \hh$. 
It is given by 
$ 
I(r)\mapsto \bar{\kappa}^{(r)}$ for all $ r\in\ZZ$ and $K_2\mapsto \bar{C}_2$,
where    $\bar{\kappa}^{(r)}$ and $\bar{C}_2$ stand for the images of Heisenberg subalgebra generators $\kappa^{(r)}$ and $C_2$  in   the corresponding component of $\text{gr}\ts \hh$.

\subsection{Iohara--Kohno realization} \label{subsec01}

We follow \cite{IK} to introduce a certain   realization    of level 1 modules for the double Yangian $\dygl$. Let $\s=\CC e_{11}\oplus\CC e_{22}$ be the Cartan subalgebra of $\gl$. Denote by $\varepsilon_j$, $j=1,2$, the linear maps $\s\to\CC$ given by $\varepsilon_j (e_{kk})=\delta_{jk}$. 
The dual $\s^*$ is equipped with the standard bilinear form $(\cdot,\cdot)$  determined by $(\varepsilon_j,\varepsilon_k)=\delta_{jk}$.
Let $Q=\ZZ\alpha$ with $\alpha=\varepsilon_1 -\varepsilon_2$ be the root lattice of the simple Lie algebra $\sll$ and $\lambda_i =\Lambda_i-\Lambda_0$ the classical part of the $i$-th fundamental affine  weight $\Lambda_i$.

Consider the   Lie algebra $\wht{\s} $    generated by the central element $C$ and the elements $a_{j}(r) $, where $j=1,2$ and $r\in\ZZ$, $r\neq 0$, subject to the relations
\beq\label{heisenberg}
\left[a_j(r),a_k(s)\right] = r\ts \delta_{jk}\ts \delta_{r+s\ts 0}\ts C.
\eeq
Let $\wht{\s}^-$ be its (commutative) subalgebra generated by all  $a_{j}(r)$ with $j=1,2$ and $r\in\ZZ_{<0}$. Denote by $\CC \left[Q\right]$    the group algebra of the root lattice $Q$.
For any $i=0,1$ and $s\in\CC$ set 
\beq\label{fis}
\mathcal{F}_{i s}=
U(\wht{\s}^-)[[h]]
\ot
\CC [Q ][[h]] e^{\lambda_i+\frac{s}{2}(\varepsilon_1 +\varepsilon_2)}.
\eeq
   The tensor product in \eqref{fis} is understood in the $h$-adically completed sense, so that the $\CC[[h]]$-module $\mathcal{F}_{i   s}$ is topologically free. Define the action of  
$C$, $a_j(r)$, $\partial_{\varepsilon_j}$ and $e^{\varepsilon_j}$, where $j=1,2$ and $r\neq 0$, on   $\mathcal{F}_{i   s}$ so that for any $f\otimes e^{\mu}\in \mathcal{F}_{i   s}$ we have
\begin{align*}
&C\cdot f\otimes e^{\mu} =f\otimes e^{\mu},\qquad\hspace{61pt}
a_j(r)\cdot f\otimes e^{\mu}=\begin{cases}
a_j(r)\ts f\otimes e^{\mu},&\text{if }r<0,\\
\left[a_j(r), f\right]\otimes e^{\mu},&\text{if }r>0,
\end{cases}\\
&\partial_{\varepsilon_j}\cdot f\otimes e^{\mu}= (\varepsilon_j,\mu)\ts f\otimes e^{\mu},\qquad\qquad
e^{\varepsilon_j}\cdot  f\otimes e^{\mu}= f\otimes e^{\varepsilon_j +\mu}.
\end{align*}
Let us  write $X_{\pm\alpha }(u)=h^{-1}X^{\pm}(u)$. 
Introduce the following operator series on $\mathcal{F}_{i   s}$:
\begin{gather*}
  E^-(u)=\exp\left(\sum_{r>0}\left( -\frac{a_1(-r)}{r}\left(u-\frac{3}{4}h\right)^r +\frac{a_2(-r)}{r}\left(u+\frac{h}{4}\right)^r\right)\right),\\
  E^+(u)=\exp\left(\sum_{r>0}\left(  \frac{a_1(r)-a_2(r)}{r} \left(u+\frac{h}{4}\right)^{-r}\right)\right),\qquad E^0(u)=e^{\alpha } \left(u+\frac{h}{4}\right)^{\partial_{\alpha }}.
\end{gather*}
Finally, we   recall the Iohara--Kohno realization \cite[Thm. 3.1]{IK}.

\begin{thm}\label{IKthm}
The following assignments define a structure of $\dygl$-module on $\mathcal{F}_{i  s}$:
\begin{align*}
X_{\alpha }(u)&\mapsto E^-(u)E^+(u)E^0(u),\\
X_{-\alpha }(u)&\mapsto  E^-(u+h/2)^{-1}E^+(u-h/2)^{-1}E^0(u-h/2)^{-1},\\
 k_{j}^{+}(u)&\mapsto \exp\left(-\sum_{r>0 }\frac{a_j(r)}{r}\left(\left(u+\frac{h}{2} \right)^{-r}- \left(u-\frac{h}{2} \right)^{-r}  \right)\right)
\left(1-\frac{h }{u+\frac{h}{2} }\right)^{\partial_{\varepsilon_j}}, j=1,2,\\
k_{j}^{-}(u)&\mapsto \exp\left(\sum_{r>0 }\left(\frac{ \delta_{1j} \ts a_2(-r)}{r}\left(\left(u+h\right)^r -u^r\right)+\frac{\delta_{2j}\ts a_1(-r)}{r}\left(u^r -(u-h)^r\right)\right)  \hspace{-4pt}\right), j=1,2.
\end{align*}
\end{thm}

\section{Semi-infinite construction}\label{sec02}

Throughout the rest of the paper, we assume that $s=0$ in \eqref{fis} and  we   write $\mathcal{F}_i=\mathcal{F}_{i\ts 0}$ for $i=0,1$. Furthermore, from now on, we assume that the double Yangian $\dys$ and the Heisenberg subalgebra $\hh$, as given in   Subsection \ref{subsec00}, are $h$-adically completed. Our goal is to construct the  topological    bases for     certain   $\dygl$-modules in parallel with the   semi-infinite construction   \cite{FS} for the affine Kac--Moody Lie algebra $\wht{\mathfrak{sl}}_2$.

\subsection{Normalized Iohara--Kohno realization} \label{subsec01}

We now introduce another structure of level $1$ $\dygl$-module over   $\mathcal{F}_{i}$ by modifying the original formulae from Theorem \ref{IKthm}.

\begin{thm}\label{IKthm2} 
The following assignments define a structure of $\dygl$-module on $\mathcal{F}_{i}$:
\begin{align*}
X_{\alpha }(u)\,&\mapsto\,  E^-(0)^{-1} E^-(u)E^+(u)E^0(u) \frac{u^{\partial_{\alpha/2} }(u+h)^{\partial_{\alpha/2}}}{(u+h/4)^{\partial_{\alpha }}},\\
X_{-\alpha }(u)\,&\mapsto \,  E^-(0)  E^-(u+h/2)^{-1}E^+(u-h/2)^{-1}E^0(u-h/2)^{-1}
 \frac{(u-h/4)^{\partial_{\alpha }}}{(u-h/2)^{\partial_{\alpha/2}}(u+h/2)^{\partial_{\alpha/2}}}
,\\
 k_{1}^{+}(u)\,&\mapsto\, \exp\left(-\sum_{r>0 }\frac{a_1(r)}{r}\left(\left(u+\frac{h}{2} \right)^{-r}- \left(u-\frac{h}{2} \right)^{-r}  \right)\right)
\left(\frac{ u+h/4}{ u+5h/4 }\right)^{\partial_{\alpha/2}}  ,\\
k_{2}^{+}(u)\,&\mapsto\, \exp\left(-\sum_{r>0 }\frac{a_2(r)}{r}\left(\left(u+\frac{h}{2} \right)^{-r}- \left(u-\frac{h}{2} \right)^{-r}  \right)\right)
\left(\frac{ u+h/4}{ u-3h/4 }\right)^{\partial_{\alpha/2}} 
\end{align*}
and the action of $k_{1}^{-}(u),k_{2}^{-}(u)$ is given by Theorem \ref{IKthm}.
\end{thm}

\begin{prf}
As with the original Iohara--Kohno realization from Theorem \ref{IKthm}, one  verifies by a direct calculation that the given   operators satisfy the defining relations for   double Yangian, which are established in  \cite[Thm. 2.1]{IK}. As a demonstration, we shall prove that the   operators satisfy the  level $1$  defining relation
\begin{align}
\left[X_\alpha (u),X_{-\alpha}(v)\right]=&\,\,\frac{1}{h}\left(
\delta(u-v-h/2)\ts k_2^+(u-h/4)\ts k_1^+(u-h/4)^{-1}\right.\non\\
&\left.-\delta(u-v+h/2)\ts k_2^-(v-h/4)\ts k_1^-(v-h/4)^{-1}
\right),\label{jps1}
\end{align}
cf. \cite[Thm. 2.1]{IK}, where $\delta(u-v)=\sum_{r\in\ZZ} u^{-r-1}v^r$ is the $\delta$-function. Throughout the proof, we  denote by $X'_{\pm\alpha}(u)$ and $k_j^{'\pm} (u)$ the corresponding operators from Theorem \ref{IKthm} in order to distinguish them from those given by Theorem \ref{IKthm2}.

First, we rewrite the left hand-side of \eqref{jps1} as
\begin{align}
\left[X_\alpha (u),X_{-\alpha}(v)\right]= 
\left[X'_\alpha (u),X'_{-\alpha}(v)\right]f(u,v),
\label{jps2}
\end{align}
where
$$
f(u,v)=\frac{u^{\partial_{\alpha/2} }(u+h)^{\partial_{\alpha/2}}}{(u+h/4)^{\partial_{\alpha }}}
\frac{(v-h/4)^{\partial_{\alpha }}}{(v-h/2)^{\partial_{\alpha/2}}(v+h/2)^{\partial_{\alpha/2}}}.
$$
Next,   consider the right hand-side of \eqref{jps1}. Using the properties of   $\delta$-function, we find
$$
\delta(u-v+h/2)
=
f(u,u+h/2)\ts\delta(u-v+h/2) =
 f(u,v)\ts \delta(u-v+h/2)  
$$
since $f(u,u+h/2)=1$. Furthermore, we have $k_j^{'-} (u)=k_j^{-} (u)$ for $j=1,2$.
Therefore, the second summand in \eqref{jps1} can be written as
\begin{align}
&-\delta(u-v+h/2)\ts k_2^-(v-h/4)\ts k_1^-(v-h/4)^{-1}\non\\
 =&-f(u,v)\ts \delta(u-v+h/2)\ts   k_2^{'-}(v-h/4)\ts k_1^{'-}(v-h/4)^{-1}.\label{jps3}
\end{align}
Finally, we consider the first summand on the right hand-side of 
\eqref{jps1}. Observe that the operators  $k_j^{+}(u)$ and $k_j^{'+}(u)$ are connected by the identities
\begin{align*}
&k_1^+(u)=
k_1^{'+}(u)\left(\frac{ u+h/2}{ u-h/2 }\right)^{\partial_{\varepsilon_1}}\left(\frac{ u+h/4}{ u+5h/4 }\right)^{\partial_{\alpha/2}},
\\
&k_2^+(u)=
k_2^{'+}(u)
\left(\frac{ u+h/2}{ u-h/2 }\right)^{\partial_{\varepsilon_2}}
\left(\frac{u+h/4}{u-3h/4}\right)^{\partial_{\alpha/2}}.
\end{align*}
Thus, we have the equality
\begin{align}
&\delta(u-v-h/2)\ts k_2^+(u-h/4)\ts k_1^+(u-h/4)^{-1}\non\\
=&\, \delta(u-v-h/2)\ts
\left(\frac{u-3h/4}{u+h/4}\right)^{\partial_{\alpha}}
\left(\frac{u+h}{u-h }\right)^{\partial_{\alpha/2}}
\ts k_2^{'+}(u-h/4)\ts k_1^{'+}(u-h/4)^{-1}.\label{jps5}
\end{align}
However, the properties of   $\delta$-function imply
\begin{align*}
&\delta(u-v-h/2) 
\left(\frac{u-3h/4}{u+h/4}\right)^{\partial_{\alpha}}
\left(\frac{u+h}{u-h }\right)^{\partial_{\alpha/2}}\\
=& \, f(u,u-h/2)\ts \delta(u-v-h/2)  =f(u,v)\ts\delta(u-v-h/2),
\end{align*}
so we conclude by \eqref{jps5} that the first summand on the right hand-side of 
\eqref{jps1} equals
\beq\label{jps4}
f(u,v)\ts
\delta(u-v-h/2)
\ts k_2^{'+}(u-h/4)\ts k_1^{'+}(u-h/4)^{-1}.
\eeq
As the operators from Theorem \ref{IKthm}   satisfy \eqref{jps1},
the commutation relation \eqref{jps1} now follows by comparing the expressions in \eqref{jps2}, \eqref{jps3} and \eqref{jps4}.
\end{prf}

The   formula for the action of $X_{\alpha }(u)$ from Theorem \ref{IKthm2} implies 
\beq\label{tzuio}
X_{\alpha }(u) \ts 1\ot 1\in 1\ot e^{\alpha}+u\mathcal{F}_{0}[[u]]
\eeq
on $\mathcal{F}_{0}$.
On the other hand, \eqref{tzuio} does not hold for the action of $X_{\alpha }(u)$ given by Theorem \ref{IKthm}. As this property is required for the semi-infinite construction in Subsection \ref{subsec15}, in the rest of the paper  we only consider   the $\dygl$-action from Theorem \ref{IKthm2} and, furthermore,  we write $X_{\pm \alpha}(u)$ and $k_j^{\pm } (u)$ for the corresponding vertex operators.

\subsection{Commutative operators} \label{subsec12}

In this subsection, we  associate with the operators $X_{  \alpha}(u)$ from Theorem \ref{IKthm2} certain commutative operators $\wvr{X} (u)$ which satisfy  an $h$-adic version of the   integrability relations which go  back to Lepowsky and Primc  \cite{LP}. As with the construction of commutative operators of Ding and B. Feigin \cite{DF}, which are associated with the Frenkel--Jing realization  \cite{FJ} of certain $U_q (\wht{\mathfrak{sl}}_2)$-modules, we   modify the term  $E^+(u)$ consisting of annihilation operators $a_j(r)$ with $r>0$. However, in contrast with the quantum affine algebra case, the term $E^0(u)$   needs to be renormalized as well.

Consider  the operator series on $\mathcal{F}_{i}$, $i=0,1$, defined by
\beq\label{xbar}
\wvr{X} (u)=X_{\alpha}(u )\ts \wvr{E}^+(u )\ts \wvr{E}^0(u ),\qquad\text{where}\qquad \wvr{E}^0(u)= (1-h/u )^{\partial_{\alpha/2} },
\eeq
the action of
$X_{  \alpha}(u)$ is given by Theorem \ref{IKthm2}
and  the term   $\wvr{E}^+(u )$ is  
\begin{align*}
&\wvr{E}^+(u)=\exp\left(\sum_{r>0}\left(\frac{a_1(r)}{r}\left(u-\frac{7}{4}h\right)^{-r} +\frac{a_2(r)}{r}\left(u+\frac{1}{4}h\right)^{-r} \right)\right)
.
\end{align*}

\begin{pro}
The following identities hold for operators on $\mathcal{F}_{i}$, $i=0,1$:
\begin{gather}
 \wvr{X} (u)\ts \wvr{X} (v)=\wvr{X} (v)\ts \wvr{X} (u),\label{comm}\\
 \wvr{X} (u)\ts \wvr{X} (u-h)=\wvr{X} (u )\ts \wvr{X} (u+h)=0.\label{int}
\end{gather}
\end{pro}

\begin{prf}
By using \eqref{heisenberg} one easily verifies the following relations:
\begin{align}
&E^+(u)E^-(v)=\frac{(u-v)(u-v+h)}{(u+h/4)^2}E^-(v)E^+(u),\label{rel1}\\
&E^0(u)E^0(v)=\frac{(u+h/4)^2}{(v+h/4)^2}E^0(v)E^0(u).\label{rel2}
\end{align}
The commutativity \eqref{comm} can be verified   using  \eqref{rel1},  \eqref{rel2}  and the identities
\begin{align}
&\wvr{E}^+(u) E^{-}(v)=\frac{u-v-h}{u-v}\frac{u+h/4}{u-7h/4}   E^{-}(v) \wvr{E}^+(u),\label{rel3}\\
&\wvr{E}^0(u)E^0(v)=\left(1-h/u\right)
  E^0(v)\wvr{E}^0(u),\label{rel4}
\end{align}
which also follow from \eqref{heisenberg}. As for the equality   \eqref{int}, it is sufficient to observe that, due to relations
\eqref{rel1}, \eqref{rel2}, \eqref{rel3} and \eqref{rel4}, we have the   decomposition 
\begin{gather*}
\wvr{X}  (u) \wvr{X}  (v)
=
f(u,v) \ts  \wvr{X}  (u,v) ,\qquad \text{where} \qquad f(u,v)=(u-v-h)(u-v+h),
\\
 \wvr{X}  (u,v) = E^-(0)^{-2} E^-(u)\ts E^-(v)\ts  E^+(u)\ts \wvr{E}^{+}(u)\ts E^+(v)\ts \wvr{E}^{+}(v)\ts  e^{2\alpha}  \ts 
  (u^2-h^2)^{\partial_{\alpha/2}}\ts (v^2-h^2)^{\partial_{\alpha/2}}
 .
\end{gather*}
By applying  $ \wvr{X}  (u,v)$ to an arbitrary element of $\mathcal{F}_{i }$ we get only finitely many negative powers of the variables $u$ and $v$ modulo $h^n$ for any   $n\geqslant 1$. Hence the formal limit $v\to u\pm h$ of  $ \wvr{X}  (u,v)$ is well-defined. Therefore,
  the limit $v\to u\pm h$ of
the product $\wvr{X}  (u) \wvr{X}  (v)$ is well-defined as well.  Furthermore, since   
$f(u,u\pm h)=0$,   it equals   zero, as required. 
\end{prf}

\subsection{Principal submodule \texorpdfstring{$W_0 $}{W0 }} \label{subsec13}

Let us write the operator series  $X_{\alpha}(u)$ and $\wvr{X} (u)$   as
$$
X_{\alpha}(u)=\sum_{r\in\ZZ} x (r)u^{-r-1}
\fand
\wvr{X} (u)=\sum_{r\in\ZZ} \wvr{x} (r)u^{-r-1}.
$$
Let $\vac=1\ot 1 \in \mathcal{F}_0$.
Motivated by the notion of principal subspaces of  affine Lie algebra modules from   \cite{FS,G}, we
define the    {\em principal submodule} $W_0 $ of  $\mathcal{F}_0$   as the $h$-adic completion of the $\CC[[h]]$-module which is spanned by all  
\beq\label{mons1}
x(r_n)\ldots x(r_1)\vac,\qquad\text{where}\qquad n\geqslant 0\fand r_1,\ldots ,r_n\in\ZZ.
\eeq
In this   subsection, we use the monomials \eqref{mons1} to construct a topological basis for $W_0$.

\begin{rem}\label{rem1}
Note that  $W_0 $ coincides with the $h$-adically  completed $\CC[[h]]$-span of all   
\beq\label{mons2}
\wvr{x}(r_n)\ldots \wvr{x}(r_1)\vac,\qquad\text{where}\qquad n\geqslant 0\fand r_1,\ldots ,r_n\in\ZZ.
\eeq
Indeed, the terms in \eqref{mons2} are   the coefficients of $u_n^{-r_n-1}\ldots u_1^{-r_1 -1}$ in $\wvr{X}(u_n)\ldots \wvr{X}(u_1)\vac$. However, by using \eqref{rel3} and \eqref{rel4}, along with $\wvr{E}^+(u)\vac =\wvr{E}^0(u)\vac=\vac$, we find
\begin{align}
\wvr{X}(u_n)\ldots \wvr{X}(u_1)\vac=F\cdot X_\alpha (u_n )\ldots X_\alpha(u_1 )\vac  \,\text{ for }\,   F=
\prod_{1\leqslant r<s\leqslant n}\left(1-\frac{h}{u_s-u_r}\right).\label{mons3}
\end{align}
As $F$ is   invertible in $\CC((u_n))\ldots ((u_1))[[h]]$ and  both sides of \eqref{mons3} belong to the $h$-adic completion of $\mathcal{F}_0((u_n))\ldots ((u_1))$, we conclude  that each element in \eqref{mons1} can be expressed  in terms  of elements of the form \eqref{mons2} and vice versa.
 Moreover, by rewriting   \eqref{mons3} as
\begin{gather*}
 h^{-1}\left(\wvr{X}(u_n)\ldots \wvr{X}(u_1)\vac-  X_\alpha (u_n )\ldots X_\alpha(u_1 )\vac\right)\\
=   h^{-1} (F-1)X_\alpha (u_n )\ldots X_\alpha(u_1 )\vac =h^{-1} (1-F^{-1})\wvr{X} (u_n)\ldots \wvr{X}(u_1)\vac,
\end{gather*}
we conclude that the expressions of the form
$ 
h^{-1}\left(x(r_n)\ldots x(r_1)\vac-\wvr{x}(r_n)\ldots \wvr{x}(r_1)\vac\right) 
$ 
can be also written as    $\CC[[h]]$-linear combinations  of elements \eqref{mons1} or \eqref{mons2}. Observe that such linear combinations can be infinite, but they are well-defined as they are finite modulo $h^m$ for all $m\geqslant 1$  and the $\CC[[h]]$-module $\mathcal{F}_0$ is $h$-adically complete.  
\end{rem}

In the following lemma, we consider the monomials of the form \eqref{mons2}. In comparison with \eqref{mons1}, they are easier to handle   since their factors $\wvr{x}(r)$ commute; recall \eqref{comm}. Let $\wvr{B}_{W_0 }$ be the set 
of
all elements \eqref{mons2} with $n\geqslant 0$ such that their indices satisfy
\beq\label{dtwo}
 r_1\leqslant -1
\Fand
r_{j+1}\leqslant r_{j} -2\quad\text{for all}\quad j=1,\ldots,n-1 .
\eeq

\begin{lem}\label{l1}
 The $h$-adic completion of the $\CC[[h]]$-span of  $\wvr{B}_{W_0 }$ coincides with $W_0 $.
\end{lem}

\begin{prf}
First, note that $\wvr{X}(u)\vac$ possesses only nonnegative powers of the variable $u$, i.e. we have $\wvr{x}(r)\vac =0 $ for $r\geqslant 0$. Thus, by commutativity \eqref{comm},   the terms 
\beq\label{span}
\wvr{x}(r_n)\ldots \wvr{x}(r_1)\vac\qquad \text{such that}\qquad r_1,\ldots ,r_n\in\ZZ_{<0}\fand r_n\leqslant\ldots \leqslant r_1.
\eeq
span an $h$-adically dense $\CC[[h]]$-submodule of $W_0 $.
Next, we  use  the $h$-adic integrability relation    to establish the difference two condition \eqref{dtwo} for the terms in \eqref{span}. Although this is just an adjustment of the well-known argument, originating from the (integrable) representation theory of $\wht{\mathfrak{sl}}_2$, to the $h$-adic setting, we give some details for completeness. Extracting the coefficients of   $u^{-2r-3}$ and $u^{-2r-2}$   in \eqref{int} and using   \eqref{comm} we find
\begin{align}
&\wvr{x}(r)\wvr{x}(r+1)=
-\sum_{l\geqslant 1}\wvr{x}(r-l)\wvr{x}(r+l+1)\mod h,\label{id1}\\
&\wvr{x}(r)\wvr{x}(r)=
-2\sum_{l\geqslant 1}\wvr{x}(r-l)\wvr{x}(r+l)\mod h,\label{id2}
\end{align}
respectively. Suppose that we have two adjacent operators $\wvr{x}(r_{j+1})\wvr{x}(r_{j})$ in the element $\wvr{x}(r_n)\ldots \wvr{x}(r_1)\vac$ of \eqref{span} which do not satisfy the difference two condition $r_{j+1}\leqslant r_{j} -2$. Then, they can be replaced modulo $h$ by a sum of those which satisfy this condition using the identity \eqref{id1} (resp. \eqref{id2}) if   $ r_{j+1}=r_j -1$  (resp. $ r_{j+1}=r_j  $).
Due to \eqref{int}, the difference of the left and the right hand-side in \eqref{id1} and in \eqref{id2} belongs to $hW_0$. In other words, it is   a $\CC[[h]]$-linear combination of monomials as in \eqref{span} such that its coefficients belong to $h\CC[[h]]$. Thus, for any   $m\geqslant 1$, by      repeating  such procedure sufficient number of times one can  express the original element $\wvr{x}(r_n)\ldots \wvr{x}(r_1)\vac$ as a $\CC[[h]]$-linear combination of elements  which satisfy the difference two condition \eqref{dtwo} modulo $h^m$. As $W_0 $ is $h$-adically complete, this implies the  lemma.
\end{prf}

In the definition \eqref{xbar} of the series  $\wvr{X} (u)$,    we multiplied    $X_{\alpha}(u )$ by the term $\wvr{E}^+(u)$, thus changing the classical limit of the resulting operator. Therefore, we can not instantly establish the linear independence of the spanning set from the previous lemma via classical limit. Instead, we now turn to the monomials of   operators \eqref{mons1}. Let $B_{W_0 }$ be the set of all elements \eqref{mons1} with $n\geqslant 0$ which satisfy the conditions imposed by \eqref{dtwo}.

\begin{lem}\label{l2}
The set $B_{W_0 }$ is linearly independent over $\CC[[h]]$.
\end{lem}

\begin{prf}
The classical limit $h\to 0$ of the 
operators $X_{\pm \alpha}(u)$ from Theorem 
\ref{IKthm2} with $i=0$
produces  the action of $\wht{\mathfrak{sl}}_2$ 
on its level 1 integrable highest weight  module $L(\Lambda_0)$ of the highest weight $\Lambda_0$, as given by the    Frenkel--Kac--Segal realization \cite{FK,S}.
Thus, the classical limit of   $B_{W_0 }$ 
produces the well-known quasi-particle basis for the principal subspace of $L(\Lambda_0)$; see \cite{FS,G}. Hence, in particular,   $B_{W_0 }$  is linearly independent over $\CC[[h]]$.
\end{prf}

Finally, we combine Lemmas \ref{l1}  and \ref{l2}  to obtain a basis of $W_0 $.

\begin{thm}\label{l3}
The set $B_{W_0 }$   forms a topological basis of $W_0 $.
\end{thm}

\begin{prf}
By Lemma \ref{l2}, it is sufficient to check that the $\CC[[h]]$-span of $B_{W_0 }$ is $h$-adically dense in $W_0 $. Choose any $w\in W_0 $. By Lemma \ref{l1}, it can be expressed as   $w=\sum_j a_j \wvr{b}_j$, where $a_j\in\CC[[h]]$ and  $\wvr{b}_j\in\wvr{B}_{W_0 }$.  Next, by Remark \ref{rem1}, we have $w=\sum_j a_j  b_j\mod h$, where $b_j\in B_{W_0 }$ denote the elements   corresponding to $\wvr{b}_j$, i.e., more precisely, if $\wvr{b}_j =\wvr{x}(r_n)\ldots \wvr{x}(r_1)\vac$, then we set $b_j =x(r_n)\ldots x(r_1)\vac$. In addition, the remark implies that $h^{-1} (w- \sum_j a_j  b_j )$ is again a $\CC[[h]]$-linear combination of the elements of $\wvr{B}_{W_0 }$. Hence, for any $m\geqslant 1$,      we can repeat this procedure until we express $w$ as a $\CC[[h]]$-linear combination of the elements of $B_{W_0 }$ modulo $h^m$, so   the theorem now follows. 
\end{prf}

\subsection{\texorpdfstring{$\dysl$}{DY(sl2)}-modules \texorpdfstring{$\mathcal{L}_i(\mathfrak{sl}_2)$}{Li(sl2)}} \label{subsec14}

Consider  the set
 $$B_{\mathcal{L}_0(\mathfrak{sl}_2)} =\left\{
e^{k\alpha} \ts b\,:\,k\in\ZZ,\,b\in B_{W_0 }  \right\} \subset   \mathcal{F}_0.$$ 
Let $\mathcal{L}_0(\mathfrak{sl}_2)$ be the $h$-adically completed $\CC[[h]]$-span of $B_{\mathcal{L}_0(\mathfrak{sl}_2)} $.  In this subsection, we show that $\mathcal{L}_0(\mathfrak{sl}_2)$ is closed with respect to the
$\dysl$-action  from Theorem \ref{IKthm2}. 

\begin{lem}\label{ll1}
For any integer $m $ we have
\beq\label{inv0}
  E^+(u)^{  m}\ts \mathcal{L}_0(\mathfrak{sl}_2)\subset \mathcal{L}_0(\mathfrak{sl}_2)[[u^{-1}]].
\eeq
\end{lem}

\begin{prf}
Suppose $m=1$.
It is sufficient to show that $E^+(u) e^{k\alpha}b$ belongs to $ \mathcal{L}_0(\mathfrak{sl}_2)[[u^{-1}]]$ for all $e^{k\alpha}b\in B_{\mathcal{L}_0(\mathfrak{sl}_2)}$.
Using the   formulae from Theorem \ref{IKthm2} and  $E^+(u)\vac =\vac$   we find
\begin{align}\label{inv1}
E^+(u)\ts e^{k\alpha}\ts X_\alpha(v_n)\ldots X_\alpha(v_1)\vac
=F\cdot  e^{k\alpha}\ts X_\alpha(v_n)\ldots X_\alpha(v_1)\vac,
\end{align} 
where, due to relation \eqref{rel1}, the factor $F\in\CC[u^{-1}][[v_n,\ldots ,v_1,h]]$ is given by
$$
F=\prod_{r=1}^n
 \left(1-\frac{ v_r }{u  }\right)   \left(1-\frac{ v_r }{u +h }\right)
 .
$$
For any integers $s_1,\ldots ,s_n$,
by extracting the coefficients
 of $v_1^{s_1}\ldots v_n^{s_n}$ from the right hand-side of \eqref{inv1}
one obtains a power series in the variable $u^{-1}$. Clearly, its coefficients are $\CC[[h]]$-linear combinations of the elements    of the form  $e^{k\alpha}x(r_n)\ldots x(r_1)\vac$. Moreover, by Theorem \ref{l3}, every   $x(r_n)\ldots x(r_1)\vac$ can be expressed in terms of elements of $B_{W_0 }$, so that every $e^{k\alpha}x(r_n)\ldots x(r_1)\vac$ belongs to $\mathcal{L}_0(\mathfrak{sl}_2)$. On the other hand, by extracting the coefficients of $v_1^{s_1}\ldots v_n^{s_n}$   from the left hand-side of \eqref{inv1} for a suitable choice of $s_1,\ldots ,s_n$, we obtain $E^+(u) e^{k\alpha}b $ with   $e^{k\alpha}b\in B_{\mathcal{L}_0(\mathfrak{sl}_2)}$. Thus, we conclude that all $E^+(u) e^{k\alpha}b$ such that $e^{k\alpha}b\in B_{\mathcal{L}_0(\mathfrak{sl}_2)}$ belong to $\mathcal{L}_0(\mathfrak{sl}_2)[[u^{-1}]]$, as required.

If $m=-1$, we can use $E^+(u)^{-1}\vac=\vac$ and verify  the inclusion in  \eqref{inv0}  analogously.
Finally, if $m$ is an arbitrary integer, one    easily proves  by induction over $|m|$  that we have
\beq\label{rwrt}
E^+(u_1)^{\pm 1}\ldots E^+(u_{|m|})^{\pm 1} \ts \mathcal{L}_0(\mathfrak{sl}_2)\subset \mathcal{L}_0(\mathfrak{sl}_2)[[u_1^{-1},\ldots ,u_{|m|}^{-1}]].
\eeq
Hence  the statement of the lemma  follows by setting $u_1=\ldots=u_{|m|}=u$ in \eqref{rwrt}.
\end{prf}

The proof of   Lemma \ref{ll1} relies on the fact that the coefficients of $E^+(u)$ are $\CC[[h]]$-linear combinations of monomials of annihilation operators $a_j(r)$ with $r>0$, which implies $E^+(u)^{\pm }\vac =\vac$. As $k_j^+(u)^{\pm 1}\vac =\vac$ for $j=1,2$, one can prove  the next lemma analogously.

\begin{lem}\label{ll3m}
For any integer $m $ and $j=1,2$ we have
$$
  k_j^+(u)^{  m}\ts \mathcal{L}_0(\mathfrak{sl}_2)\subset \mathcal{L}_0(\mathfrak{sl}_2)[[u^{- 1}]].
$$
\end{lem}
 
By combining Lemma \ref{ll3m} and the identity \eqref{hplus} we get

\begin{kor}\label{kor1}
$H^+(u)\ts \mathcal{L}_0(\mathfrak{sl}_2) \subset \mathcal{L}_0(\mathfrak{sl}_2) [[u^{-1}]].$
\end{kor}

Let us consider  the remaining   operators from Theorem \ref{IKthm2}  whose action employs the creation operators $a_j(r)$ with $r<0$ as well.
First,  observe that 
\beq\label{inv7}
X_\alpha(u)\ts e^{k\alpha}=u^{2k}(u+h )^{2k}\ts e^{k\alpha}\ts X_\alpha(u)
\qquad\text{implies}\qquad
X_\alpha(u)\ts \mathcal{L}_0(\mathfrak{sl}_2)\subset \mathcal{L}_0(\mathfrak{sl}_2)[[u^{\pm 1}]].
\eeq
Next, write $\mathcal{E}^{-}(u)=E^-(0)^{-1}E^{-}(u)$, so that, by Theorem \ref{IKthm2},   we have
\beq\label{cmp}
X_\alpha(u)=\mathcal{E}^{-}(u)\ts E^+(u)\ts e^\alpha\ts  u^{\partial_{\alpha/2}}(u+h)^{\partial_{\alpha/2}}.
\eeq

\begin{lem}\label{ll2}
For any integer $m $ we have
\beq\label{inv4}
  \mathcal{E}^{-}(u)^{  m}\ts \mathcal{L}_0(\mathfrak{sl}_2)\subset \mathcal{L}_0(\mathfrak{sl}_2)[[u]].
\eeq
\end{lem}

\begin{prf}
Let $m=1$.
We can use the formulae from Theorem \ref{IKthm2} and \eqref{rel1} to prove  
\begin{align}\label{inv2}
\mathcal{E}^{-}(u)\ts e^{k\alpha}\ts X_\alpha(v_n)\ldots X_\alpha(v_1)\vac
=F\cdot  e^{k\alpha}\ts X_\alpha(v_n)\ldots X_\alpha(v_1)\ts \mathcal{E}^{-}(u)\vac,
\end{align} 
where the factor $F\in \CC[v_n^{-1},\ldots, v_1^{-1}][[u,h]]$ is  given by
$$
F=\prod_{r=1}^n  \left(1-\frac{u}{v_r}\right)^{-1}\left(1-\frac{u}{v_r+h}\right)^{-1}  .
$$
Since $E^+(u)e^\alpha u^{\partial_{\alpha/2}}(u+h)^{\partial_{\alpha/2}}\vac=e^{\alpha}\vac$, by \eqref{cmp} the rightmost term in \eqref{inv2} equals
$$
\mathcal{E}^-(u)\vac=e^{-\alpha}\ts  \mathcal{E}^-(u)\ts E^+(u)\ts e^\alpha \ts u^{\partial_{\alpha/2}}(u+h)^{\partial_{\alpha/2}}\vac=e^{-\alpha}\ts  X_\alpha(u)\vac.
$$
Therefore, the right hand-side of \eqref{inv2} is equal to
\beq\label{inv3}
F\cdot  e^{k\alpha}\ts X_\alpha(v_n)\ldots X_\alpha(v_1)\ts e^{-\alpha} X_\alpha(u)\vac
=
G\cdot  e^{(k-1)\alpha}\ts X_\alpha(v_n)\ldots X_\alpha(v_1)\ts   X_\alpha(u)\vac.
\eeq
The right hand-side of \eqref{inv3} was found by moving the operator $e^{-\alpha}$ all the way to the left, which produced  the term $G\in\CC[v_n^{-1},\ldots, v_1^{-1}][[u,h]]$  given by
$$
G=F\cdot \prod_{r=1}^n  v_r^{-1}\left(v_r+h \right)^{-1}=\prod_{r=1}^n  \left(v_r-u\right)^{-1}\left(v_r -u+h\right)^{-1} .
$$
Thus, we proved that the left hand-side in \eqref{inv2} coincides with the right hand-side in \eqref{inv3}. 
However, it is clear that all  coefficients of the right hand-side of \eqref{inv3} belong to $\mathcal{L}_0(\mathfrak{sl}_2)$.
Hence,  the inclusion \eqref{inv4} for $m=1$ follows as in the proof of Lemma \ref{ll1}. Furthermore,    the general case $m>1$    can be again verified by induction over $m$. 

Let $m=-1$. Note that $\mathcal{E}^-(u)$ can be written as
$
\mathcal{E}^{-}(u) =1-e^{-}(u)$, where $e^{-}(u)$ belongs to $u U(\wht{\s}^-)[[h,u]]$.
Hence $\mathcal{E}^{-}(u)$ is invertible and its inverse  takes the form
\beq\label{cmp2}
\mathcal{E}^{-}(u)^{-1} =\left(1-e^{-}(u)\right)^{-1}=\sum_{l\geqslant 0}e^{-}(u)^l=\sum_{l\geqslant 0}\left(1-\mathcal{E}^{-}(u)\right)^l.
\eeq
Each summand in \eqref{cmp2} is given in terms of nonnegative powers of $\mathcal{E}^{-}(u)$, which satisfy \eqref{inv4}, so we conclude that  $E^{-}(u)^{-1}\mathcal{L}_0(\mathfrak{sl}_2)\subset \mathcal{L}_0(\mathfrak{sl}_2) [[u]]$, as required. Finally, as before, this is     generalized to any negative integer $m$  by induction.
\end{prf}
 
By using \eqref{hplus} and   the action of $k_j^{-}(u)$ from Theorem  \ref{IKthm2} one easily verifies 
 the   identity 
$
H^-(u-h/4) =\mathcal{E}^-(u)\ts \mathcal{E}^-(u+h)^{-1}
$ for operators on $\mathcal{F}_0$.
Hence, by  Lemma \ref{ll2} we have
\begin{kor}\label{kor2}
$H^-(u)\ts \mathcal{L}_0(\mathfrak{sl}_2) \subset \mathcal{L}_0(\mathfrak{sl}_2) [[u]].$
\end{kor}

Let us  turn our attention to the   operator  $X_{-\alpha}(u)$, which takes the form
\beq\label{cmp3}
X_{-\alpha}(u)=
\mathcal{E}^-(u+h/2)^{-1}\ts E^+(u-h/2)^{-1}\ts  e^{-\alpha}\ts 
(u-h/2)^{-\partial_{\alpha/2}}(u+h/2)^{-\partial_{\alpha/2 }}.
\eeq

\begin{lem}\label{ll3}
$X_{-\alpha}(u)\ts \mathcal{L}_0(\mathfrak{sl}_2) \subset \mathcal{L}_0(\mathfrak{sl}_2) [[u^{\pm 1}]].$
\end{lem}

\begin{prf}
The lemma follows by an argument which goes in parallel with the proof of Lemma \ref{ll1}. More specifically, it relies on the identity 
\begin{align}\label{inv6}
X_{-\alpha}(u)\ts e^{k\alpha}\ts X_\alpha(v_n)\ldots X_\alpha(v_1)\vac
=F\cdot \mathcal{E}^{-}(u+h/2)^{-1}\ts  e^{(k-1)\alpha}\ts X_\alpha(v_n)\ldots X_\alpha(v_1)\vac,
\end{align} 
where the series $F\in\CC[u^{-1}][[v_n ,\ldots ,v_1,h]]$ is given by
$$
F=\left(u-h/2\right)^{-k}\left(u+h/2\right)^{-k}\prod_{r=1}^n\left(u-v_r-h/2\right)^{-1}\left(u-v_r+h/2\right)^{-1}.
$$
In addition, it employs Lemma \ref{ll2}, which implies that all coefficients of the right hand-side of \eqref{inv6} belong to $\mathcal{L}_0(\mathfrak{sl}_2)$. 
As for the equality in \eqref{inv6}, it is   verified using the  expressions \eqref{cmp} and \eqref{cmp3} for the operators $X_{\pm\alpha}(u)$ and the relation \eqref{rel1}.
\end{prf}

Finally,   Corollaries \ref{kor1} and \ref{kor2}, the inclusion in \eqref{inv7} and Lemma  \ref{ll3}    imply

\begin{thm}
The $\CC[[h]]$-module $\mathcal{L}_0(\mathfrak{sl}_2)$ is a module for the double Yangian $\dysl$.
\end{thm}

Consider the $\CC[[h]]$-module  
$
\mathcal{L}_1(\mathfrak{sl}_2)=\left\{e^{\lambda_1}v\,:\, v\in \mathcal{L}_0(\mathfrak{sl}_2)\right\}
$.
The next lemma  is verified by a direct calculation.

\begin{lem}\label{cmr1}
The following commutation relations hold:
\begin{align*}
&X_\alpha(u)\ts e^{\lambda_1}= u^{1/2}(u+h)^{1/2}\ts  e^{\lambda_1}\ts X_\alpha(u),\\
&X_{-\alpha}(u)\ts e^{\lambda_1}= (u-h/2)^{-1/2}(u+h/2)^{-1/2}\ts   e^{\lambda_1}  X_{-\alpha}(u),\\
& H^+(u)\ts e^{\lambda_1}= (u+7h/4)^{1/2}   (u-h/4)^{-1/2}e^{\lambda_1}\ts  H^+(u), \quad
 H^-(u)\ts e^{\lambda_1}= e^{\lambda_1}\ts  H^-(u).
\end{align*}
\end{lem}

By using Lemma \ref{cmr1}, one easily checks that 
$\mathcal{L}_1(\mathfrak{sl}_2)$ is closed under the action of the double Yangian $\dysl$ as well, so that  we have 

\begin{thm}
The $\CC[[h]]$-module $\mathcal{L}_1(\mathfrak{sl}_2)$ is a module for the double Yangian $\dysl$.
\end{thm}

\subsection{Semi-infinite monomial bases for \texorpdfstring{$\mathcal{L}_i(\mathfrak{sl}_2)$}{Li(sl2)}} \label{subsec15}
In this subsection, we construct   topological bases for    $\mathcal{L}_0(\mathfrak{sl}_2)$ and $\mathcal{L}_1(\mathfrak{sl}_2)$. Let us start with $\mathcal{L}_0(\mathfrak{sl}_2)$. Introduce the operator
\beq\label{cmp6}
\wtld{X}(u)=\sum_{r\in\ZZ}\wtld{x}(r) u^{-r-1}=X_\alpha(u)\ts \left(1+h/u\right)^{-\partial_{\alpha/2}},
\eeq
where the action of $X_\alpha(u)$ is given by Theorem \ref{IKthm2}. It satisfies the identity
\beq\label{cmp4}
\wtld{X}(u_n)\ldots\wtld{X}(u_1)\vac
=
F\cdot X_\alpha(u_n)\ldots X_\alpha(u_1)\vac
\quad\text{for}\quad
F=\prod_{r=2,\ldots ,n} \left(1+ h/u_r\right)^{1-r}.
\eeq
Let $\wtld{B}_{W_0 }$ be the set 
of
all elements 
$
\wtld{x}(r_n)\ldots \wtld{x}(r_1)\vac$,  where $n\geqslant 0$ and the  integers $r_1,\ldots ,r_n$
  satisfy the   conditions given by \eqref{dtwo}.

\begin{kor}\label{cmp5}
The set $\wtld{B}_{W_0 }$   forms a topological basis of $W_0 $.
\end{kor}

\begin{prf}
The given set is linearly independent as its classical limit coincides with the classical limit of the basis $B_{W_0 }$ from Theorem \ref{l3}. The fact that the $\CC[[h]]$-span of $\wtld{B}_{W_0 }$ is $h$-adically dense in $W_0 $ is established by arguing as in the proof of Theorem \ref{l3}. However, while the argument therein relies on
 the identity \eqref{mons3},  here one   uses \eqref{cmp4} instead.
\end{prf}

For any integer $m$ let $\wtld{B}_{\mathcal{L}_0(\mathfrak{sl}_2),m}   $ be the set of all elements $e^{m\alpha}b$ such that  $b\in  \wtld{B}_{W_0 }$.
   Clearly,  the union $\cup_{m\in\ZZ} \wtld{B}_{\mathcal{L}_0(\mathfrak{sl}_2),m} $ spans an $h$-adically dense $\CC[[h]]$-submodule of $\mathcal{L}_0(\mathfrak{sl}_2)$.    We 
	now  employ the Feigin--Stoyanovsky-type 
	construction \cite{FS} to reduce it   to a linearly independent set.
By \eqref{tzuio} and \eqref{cmp6} we have
\beq\label{tld7}
\wtld{x}(-1)\vac=e^{\alpha}\vac.
\eeq
Furthermore, by extracting the coefficients in the   relation
$
   \wtld{X}(u)\ts e^{-\alpha}= e^{-\alpha}\ts \wtld{X}(u) \ts u^{-2},  
$
which can be easily verified by a direct calculation,
we find
\beq\label{tld7u}
 \wtld{x}(r)\ts e^{-\alpha} =e^{-\alpha}\ts  \wtld{x}(r-2 )   \quad\text{for all}\quad r\in \ZZ.
\eeq
Let $\mathcal{L}_0(\mathfrak{sl}_2)_m$ be the $h$-adically completed $\CC[[h]]$-span of  $\wtld{B}_{\mathcal{L}_0(\mathfrak{sl}_2),m} $. Observe that the union $\cup_{m\in\ZZ} \mathcal{L}_0(\mathfrak{sl}_2)_m $ coincides with $\mathcal{L}_0(\mathfrak{sl}_2)$. Also, we have $\mathcal{L}_0(\mathfrak{sl}_2)_m\subset \mathcal{L}_0(\mathfrak{sl}_2)_{m-1}$ for all integers $m$. Indeed,   this is an immediate consequence of the inclusion $\wtld{B}_{\mathcal{L}_0(\mathfrak{sl}_2),m} \subset \wtld{B}_{\mathcal{L}_0(\mathfrak{sl}_2),m-1}$, which can be proved by using the identities  \eqref{tld7} and \eqref{tld7u} as follows: 
\begin{align*}
e^{m\alpha}\ts \wtld{x}(r_n)\ldots  \wtld{x}(r_1)\vac
=&\, \,e^{m\alpha}\ts \wtld{x}(r_n)\ldots  \wtld{x}(r_1)\ts e^{-\alpha}\ts e^{\alpha}\vac 
=  e^{m\alpha}\ts \wtld{x}(r_n)\ldots  \wtld{x}(r_1)\ts e^{-\alpha}\ts \wtld{x}(-1)\vac\\
= & \,\,e^{(m-1)\alpha}\ts \wtld{x}(r_n -2)\ldots  \wtld{x}(r_1-2) \ts  \wtld{x}(-1)\vac \in  \wtld{B}_{\mathcal{L}_0(\mathfrak{sl}_2),m-1} .
\end{align*}
Finally, we obtain the topological basis for $\mathcal{L}_0(\mathfrak{sl}_2)$.
\begin{thm}\label{cmpdl1}
The direct limit  $$\wtld{B}_{\mathcal{L}_0(\mathfrak{sl}_2)}=\lim_{\longrightarrow}  \wtld{B}_{\mathcal{L}_0(\mathfrak{sl}_2),m}    $$  forms a topological basis of $\mathcal{L}_0(\mathfrak{sl}_2)$.
\end{thm}

 \begin{prf}
It is clear that the elements of $\wtld{B}_{\mathcal{L}_0(\mathfrak{sl}_2)}$ span  the $h$-adically dense $\CC[[h]]$-submodule of $\mathcal{L}_0(\mathfrak{sl}_2)$. On the other hand, their classical limits produce the semi-infinite basis for the integrable highest weight $\wht{\mathfrak{sl}}_2$-module $L(\Lambda_0)$ of the highest weight $\Lambda_0$,   established by B. Feigin and Stoyanovsky  \cite{FS}, so that they are linearly independent.
\end{prf}

Consider the $\dysl$-module   $\mathcal{L}_1(\mathfrak{sl}_2)$.
For any integer $m$ let $\wtld{B}_{\mathcal{L}_1(\mathfrak{sl}_2),m}   $ be the set of all elements $e^{\lambda_1+m\alpha }b$ with $b\in  \wtld{B}_{W_0 }$.
The union $\cup_{m\in\ZZ} \wtld{B}_{\mathcal{L}_1(\mathfrak{sl}_2),m} $  spans  an $h$-adically dense $\CC[[h]]$-submodule of $\mathcal{L}_1(\mathfrak{sl}_2)$.
The   semi-infinite basis for $\mathcal{L}_1(\mathfrak{sl}_2)$ is established in parallel with the case of $\mathcal{L}_0(\mathfrak{sl}_2)$. Naturally, in this case, its linear independence  follows from the observation that its classical limit  produces the semi-infinite basis for the  integrable highest weight $\wht{\mathfrak{sl}}_2$-module $L(\Lambda_1)$ of the highest weight $\Lambda_1$ from \cite{FS}. Hence, we have

\begin{thm}\label{cmpdl2}
The direct limit  $$\wtld{B}_{\mathcal{L}_1(\mathfrak{sl}_2)}=\lim_{\longrightarrow}  \wtld{B}_{\mathcal{L}_1(\mathfrak{sl}_2),m}    $$  forms a topological basis of $\mathcal{L}_1(\mathfrak{sl}_2)$.
\end{thm}

\begin{rem}
In this remark, we  discuss an interpretation of the bases $\wtld{B}_{\mathcal{L}_i(\mathfrak{sl}_2)}$ in terms of   semi-infinite monomials    \cite{FS}.
First, let us write
$v_{i,m} =e^{\lambda_i+m\alpha}\vac$ for $m\in\ZZ$ and $i=0,1$ so that, in particular, we have $v_{i,0}  =e^{\lambda_i }\vac$. We can express the elements $v_{i,0} =e^{\lambda_i }\vac$ as
$$
 v_{i,0} 
=e^{\lambda_i }\vac
=e^{\lambda_i  -\alpha}e^{\alpha}\vac
=e^{\lambda_i  -\alpha}\wtld{x}(-1)\vac
=\wtld{x}(1-i) e^{\lambda_i  -\alpha}\vac
=\wtld{x}(1-i) v_{i,-1}.
$$
Next, by a similar calculation,  which starts with $v_{i,-1}$, we get
$$
   v_{i,-1} 
	=e^{\lambda_i  -\alpha}\vac
=  e^{\lambda_i  -2\alpha}e^{\alpha}\vac
=  e^{\lambda_i -2\alpha}\wtld{x}(-1)\vac
=  \wtld{x}(3-i)e^{\lambda_i  -2\alpha}\vac
=\wtld{x}(3-i)v_{i,-2} .
$$
Combining the above equalities we find that $v_{i,0} = \wtld{x}(1-i) \wtld{x}(3-i)v_{i,-2}  $.
By repeating such
 calculations we find
\beq\label{cmpdl3} 
v_{i,0} 
=
\wtld{x}(1-i) \wtld{x}(3-i)\ldots \wtld{x}(2m-1-i)\ts v_{i,-m} \quad \text{for any}\quad  m>0.
\eeq
Note that by \eqref{cmpdl3}   the direct limits in Theorems \ref{cmpdl1} and \ref{cmpdl2} correspond  with taking the limit   $m\to \infty$ of the elements of the form
$$
e^{ m\alpha}\ts b\ts v_{i,0} =e^{ m\alpha}\ts b\ts \wtld{x}(1-i) \wtld{x}(3-i)\ldots \wtld{x}(2m-1-i)\ts v_{i,-m} ,\quad\text{where}\quad b\in  \wtld{B}_{W_0 }
$$
Hence, arguing as above to move the term $e^{ m\alpha}$ to the right, we see  that the elements of the bases $\wtld{B}_{\mathcal{L}_i(\mathfrak{sl}_2)}$  can be represented as the {\em semi-infinite monomials} 
$$
\wtld{x}(r_1)\wtld{x}(r_2)\ldots  v_{i,-\infty }\qquad\text{such that}
$$
\begin{enumerate}
\item Their   degrees $r_1,r_2,\ldots \in\ZZ$  satisfy the difference two condition
$$
r_j\leqslant r_{j+1}-2\quad\text{for all}\quad j=1,2,\ldots ,
$$
\item For each monomial there exists an index $n$ such that all degrees 
$ r_n,r_{n+1},\ldots  $ are consecutive odd (resp. even) integers if $i=0$ (resp. $i=1$).
\end{enumerate}
\end{rem}

\subsection{Semi-infinite monomial bases for \texorpdfstring{$\mathcal{L}_i(\mathfrak{gl}_2)$}{Li(gl2)}} \label{subsec31}

Consider the   $\dyg$-modules
$$
\mathcal{L}_0(\mathfrak{gl}_2)= \dyg\cdot \vac\fand \mathcal{L}_1(\mathfrak{gl}_2)= \dyg\cdot e^{\lambda_1} ,
$$
where, as before,  the action of the double Yangian is given by Theorem \ref{IKthm2}. Recall that the Heisenberg subalgebra $\hh$ commutes with all elements of $\dys$, so that we have
$$
\mathcal{L}_0(\mathfrak{gl}_2)= \hh\cdot\dys\cdot \vac =\dys\cdot\hh\cdot \vac \fand \mathcal{L}_1(\mathfrak{gl}_2)= \hh\cdot\dys\cdot e^{\lambda_1}=\dys\cdot\hh\cdot e^{\lambda_1}.
$$

Consider the subalgebra $\hh^-\subset \hh$ generated by $1$ and all elements $\kappa^{(-r)}$ with $r\geqslant 1$.
The algebra $\hh^-$ possesses the Poincar\'{e}--Birkhoff--Witt  basis
\beq\label{bhminus}
B_{\hh^-}
=
\left\{
\kappa^{(-r_n)}\ldots \kappa^{(-r_1)}\,:\, n\geqslant 0,\, r_n\geqslant\ldots \geqslant r_1\geqslant 1
\right\}.
\eeq
Indeed, 
the elements of $B_{\hh^-}$ span $\hh^-$  as all  $\kappa^{(-r)}$ with $r\geqslant 1$ commute. On the other hand, their linear independence is established by the map \eqref{hizo} and the Poincar\'{e}--Birkhoff--Witt theorem for $U(\wht{\h}^-)$, where $\wht{\h}^-$ stands for the Lie subalgebra of $\wht{\h}$ generated by all $I(-r)=I\ot t^{-r}$ with $r\geqslant 1$.
 For $i=0,1$ introduce the sets
$$
\wtld{B}_{\mathcal{L}_i(\mathfrak{gl}_2)}=\left\{h\ts b\,:\, h\in B_{\hh^-},\, b\in \wtld{B}_{\mathcal{L}_i(\mathfrak{sl}_2)}\right\}.
$$

\begin{thm}
Let $i=0,1$.
 The set $\wtld{B}_{\mathcal{L}_i(\mathfrak{gl}_2)}$ forms a topological basis for $\mathcal{L}_i(\mathfrak{gl}_2)$.
\end{thm}

\begin{prf}
Theorems \ref{cmpdl1} and \ref{cmpdl2}, along with the discussion preceding this theorem, imply that the set $\wtld{B}_{\mathcal{L}_i(\mathfrak{gl}_2)}$ spans an $h$-adically dense $\CC[[h]]$-submodule of $\mathcal{L}_i(\mathfrak{gl}_2)$. On the other hand,  its linear independence   follows   from the observation that the   classical limit of
	$\mathcal{L}_i(\mathfrak{gl}_2)$ is equal to
	$U(\wht{\h}^-)\ot L(\Lambda_i)$.
\end{prf}

\section{Modules for the quantum  vertex algebra \texorpdfstring{$\mathcal{V}^c(\mathfrak{gl}_2)$}{Vc(gl2)}}\label{sec03}

In this section, we study the underlying (quantum) vertex algebraic framework of the semi-infinite   construction. In particular, we  obtain examples of irreducible modules for the quantized universal affine vertex algebra of $\mathfrak{gl}_2$.

\subsection{Quantum affine vertex algebra \texorpdfstring{$\mathcal{V}^c(\mathfrak{gl}_2)$}{Vc(gl2)}} \label{subsec02}
The {\em dual Yangian} $\text{Y}^+(\mathfrak{gl}_2)$ for $\mathfrak{gl}_2$ is defined as the    associative algebra over the   ring $\CC[[h]]$ generated by the elements $t_{ij}^{(- r)}$, where $i,j=1,2$ and $r=1,2\ldots.$ Its defining relations are given by
\begin{align*}
R_{12}(u-v)\ts T_{13}^{-}(u)\ts T_{23}^{-}(v) &
=T_{23}^{-}(v)\ts T_{13}^{-}(u)\ts R_{12}(u-v), 
\end{align*}
where   the notation  is as in \eqref{tu} and \eqref{tij}. We use the same symbols for the generators of the dual Yangian and the corresponding generators of the double Yangian as $\text{Y}^+(\mathfrak{gl}_2)$ can be naturally regarded as a subalgebra of $\dyg$ due to the  Poincar\'{e}--Birkhoff--Witt theorem.
Moreover, for any $c\in\CC$ the $h$-adic completion of the dual Yangian is naturally equipped with the $\dyg$-action so that the central element $C$ acts as the scalar multiplication by $c$.
We shall denote this $\dyg$-module by $\mathcal{V}^c(\mathfrak{gl}_2)$.
Furthermore, $\mathcal{V}^c(\mathfrak{gl}_2)$ can be   equipped with the  quantum vertex algebra structure via  Etingof--Kazhdan's  construction  as follows; see \cite[Thm. 2.3]{EK} for more information. 

\begin{thm}\label{qvathm}
There exists a unique   quantum vertex algebra structure on $\mathcal{V}^c(\mathfrak{gl}_2)$ so that the unit $1$ is the vacuum vector and the vertex operator map satisfies
\begin{align*}
Y(T^-_{1n}(u_1)\ldots T^-_{n-1\ts n}(u_{n-1}),z)=
&\,\,T^-_{1n}(z+u_1)\ldots T^-_{n-1\ts n}(z+u_{n-1})\\
&\times T^+_{n-1\ts n}(z+u_{n-1}-hc/2)^{-1}\ldots T^+_{1  n}(z+u_{1}-hc/2)^{-1}.
\end{align*}
\end{thm}

\begin{prf}
In comparison with the original result \cite[Thm. 2.3]{EK}, we use  a different normalization of the Yang $R$-matrix, which governs the relation \eqref{rtt2} between the operators $T^{\pm}(u)$. However, the  theorem can be again proved by the analogous arguments, which can be also recovered from the proofs of \cite[Thm. 2.3.8]{G} and \cite[Thm. 4.1]{JKMY}. 
\end{prf}

\begin{rem}
The   map $Y(\cdot, z)$ can be  expressed modulo $h$ as follows. Write 
$$t^{\pm }(u)=\pm h^{-1} (I-T^\pm (u))\fand t(u)=t^{-}(u)+t^+(u).$$
One can prove by induction over $n$ the identity
$$
Y(t^-_{1n}(u_1)\ldots t^-_{n-1\ts n}(u_{n-1}),z)= \norm{t_{1n}(u_1)\ldots t_{n-1\ts n}(u_{n-1})}\mod h,
$$
where  the normal-ordered product of operators is defined   in a usual way:
$$
\norm{t(u_1)t(u_2)} =t^-(u_1)t(u_2)+t(u_2)t^+(u_1)\fand\hspace{-3pt}
\norm{t(u_1)\ldots t(u_r)} 
= \norm{ t(u_1)\hspace{-4pt}\norm{t(u_2)\ldots t(u_r)}\hspace{-2pt}}\hspace{-5pt}.
$$
Thus, it is easy to see that  the classical limit $h\to 0$ of the quantum vertex algebra $\mathcal{V}^c(\mathfrak{gl}_2)$ coincides with  the universal affine vertex algebra $V^c(\mathfrak{gl}_2)$ which goes back to the papers of I. Frenkel and Zhu \cite{FZ} and Lian \cite{Lian}; recall Remark \ref{klasicnilimes}.
\end{rem}

Let $c_2=-2c$. Denote by $\mathcal{V}^{c_2}(\mathfrak{h} )$ the
 level $c_2$
  module for  the Heisenberg subalgebra $\hh\subset \dygl$ defined over the $h$-adic completion of $\hh^-$.
	Clearly, we have
	$\mathcal{V}^{c_2}(\mathfrak{h} )\subset \mathcal{V}^c(\mathfrak{gl}_2)$ and the topological basis of $\mathcal{V}^{c_2}(\mathfrak{h} )$ is given by \eqref{bhminus}.
	The  generator series $K^\pm (u)$ of $\hh$ can be expressed as the quantum determinants of the   matrices $T^{\pm} (u)$ from the $RTT$-realization of $\dygl$, so that we have $K^\pm (u) =\qdet T^{\pm} (u)$; see \cite[Thm. B.15]{I}.
Using this observation, one obtains the following simple consequence of Theorem \ref{qvathm}.

\begin{kor} \label{korhei}
The $\CC[[h]]$-module $\mathcal{V}^{c_2}(\mathfrak{h} )$ is a  quantum vertex subalgebra of $\mathcal{V}^c(\mathfrak{gl}_2)$ for $c=-c_2/2$. Its vacuum vector is $1\in \mathcal{V}^c(\mathfrak{gl}_2)$   and its vertex operator map satisfies
\begin{align*}
Y(K^-(u_1)\ldots K^-(u_{n}),z)=
&\,\,K^-(z+u_1)\ldots K^- (z+u_{n})\\
&\times K^+ (z+u_{n}+hc_2/4)^{-1}\ldots K^+(z+u_{1}+hc_2/4)^{-1}.
\end{align*}
\end{kor}

\begin{rem}
The results of  this subsection are presented for  $\mathfrak{gl}_2$ in order to better fit the setting of the paper. However, their generalization to  $\mathfrak{gl}_N$ with $N> 2$ is straightforward. This is also true for the notion of restricted module and Theorem \ref{thm33} which we give  below.
\end{rem}

\subsection{Restricted modules for the double Yangian} \label{subsec32}
A $\dyg$-module $V$ is said to be {\em restricted} if it is a topologically free $\CC[[h]]$-module such that for any $v\in V$ and $n\geqslant 1$ the expression $T^+(u)v$ possesses only finitely many negative powers of $u$ modulo $h^n$.  The last requirement can be equivalently expressed as
\beq\label{restricted}
T^+(u)v \in\ndo\CC^2 \ot V[u^{-1}]_h\quad\text{for all}\quad v\in V,
\eeq
where $V[u^{-1}]_h$ stands for the $h$-adic completion of the $\CC[[h]]$-module $V[u^{-1}]$. Extending this notation, we shall also write $V((u))_h$ for the $h$-adic completion of $V((u))$.

\begin{pro}\label{prop32}
The $\dyg$-modules $\mathcal{F}_{i  s}$ and $\mathcal{F}_i$ established by Theorems \ref{IKthm} and \ref{IKthm2} respectively are restricted.
\end{pro}

\begin{prf}
It is clear that the underlying $\CC[[h]]$-module structure is topologically free. Moreover, by examining the explicit formulae in the aforementioned theorems we find
$$
 X^{\pm }(u)\in\om (F,F((u))_h), \quad
k_j^{-}(u)\in\om (F,F[[u]]), \quad
k_j^{+}(u)\in\om (F,F[u^{-1}]_h),
$$
 where $F=\mathcal{F}_{i  s},\mathcal{F}_i$. We now employ the identities in \eqref{rels1}--\eqref{rels3} to prove the proposition. First, we note that by the first equality in \eqref{rels1} we have 
$t_{11}^+(u)\in\om (F,F[u^{-1}]_h)$. Next, from \eqref{rels3} we obtain
$$
t_{21}^+(u)=
X^-(u-hC/4)\ts t_{11}^+(u )+t_{21}^-(u-hC/2)\ts t_{11}^- (u-hC/2)^{-1}\ts t_{11}^+(u ),  
$$
which implies $t_{21}^+(u)\in\om (F,F[u^{-1}]_h)$. As for   
$t_{12}^+(u)$, we rewrite the relation \eqref{rels2} as
$$
t_{11}^+(u)^{-1}\ts t_{12}^+(u) = X^+(u+hC/4)
+t_{11}^-(u+hC/2)^{-1}\ts t_{12}^-(u+hC/2)
$$
to conclude that the product $t_{11}^+(u)^{-1}  t_{12}^+(u)$ belongs to
$\om (F,F[u^{-1}]_h)$. Multiplying this product by $t_{11}^+(u)\in \om (F,F[u^{-1}]_h)$, we obtain $t_{12}^+(u)\in \om (F,F[u^{-1}]_h)$, as required. Finally, the second  equality in \eqref{rels1} implies
$$
t_{22}^{+}(u)= k_2^{+}(u)+ t_{21}^{+}(u)\ts t_{11}^{+}(u)^{-1}\ts t_{12}^{+} (u).
$$
As $t_{11}^{+}(u)^{-1}\in \om (F,F[u^{-1}]_h)$, we conclude that the remaining matrix entry $t_{22}^{+}(u)$ of $T^+(u)$ belongs to 
$\om (F,F[u^{-1}]_h)$, thus verifying the requirement imposed by \eqref{restricted}.
\end{prf}

The next theorem was our main motivation for introducing the notion of restricted module. Its proof  relies on the $RTT$-realization of the double Yangian and goes in parallel with the proof of Theorem \ref{qvathm}.

\begin{thm}\label{thm33}
Let $V$ be a restricted $\dyg$-module of level $c$. Then there exists a unique structure of $\mathcal{V}^c(\mathfrak{gl}_2)$-module over $V$ such that for all $n\geqslant 1$ we have
\begin{align*}
Y_V(T^-_{1n}(u_1)\ldots T^-_{n-1\ts n}(u_{n-1}),z)=
&\,\,T^-_{1n}(z+u_1)_V\ldots T^-_{n-1\ts n}(z+u_{n-1})_V\\
&\times T^+_{n-1\ts n}(z+u_{n-1}-hc/2)_V^{-1}\ldots T^+_{1  n}(z+u_{1}-hc/2)_V^{-1}.
\end{align*}
\end{thm}

The next corollary is an immediate consequence of Proposition \ref{prop32} and Theorem \ref{thm33}.

\begin{kor}\label{thm33cor}
There exists a unique structure of $\mathcal{V}^1(\mathfrak{gl}_2)$-module over $\mathcal{L}_i= \mathcal{L}_i(\mathfrak{gl}_2)$ with $i=0,1$  such that for all $n\geqslant 1$ we have
\begin{align*}
Y_{\mathcal{L}_i}(T^-_{1n}(u_1)\ldots T^-_{n-1\ts n}(u_{n-1}),z)=
&\,\,T^-_{1n}(z+u_1)_{\mathcal{L}_i}\ldots T^-_{n-1\ts n}(z+u_{n-1})_{\mathcal{L}_i}\\
&\times T^+_{n-1\ts n}(z+u_{n-1}-h/2)_{\mathcal{L}_i}^{-1}\ldots T^+_{1  n}(z+u_{1}-h/2)_{\mathcal{L}_i}^{-1}.
\end{align*}
\end{kor}

Consider the standard structure of bosonic Fock    $\wht{\h}$-module  $M(-2)$ of level $-2$. In particular, recall that $M(-2)$ is an irreducible $\wht{\h}$-module which, as a  vector space, coincides with $U(\wht{\h}^{-})$. The tensor product 
$ L_i =M(-2)\ot L(\Lambda_i)$ for $i=0,1 $,  
where $L(\Lambda_i)$  
is  the  integrable highest  weight  $\wht{\mathfrak{sl}}_2$-module  of highest weight  $\Lambda_i$,
 is naturally equipped with the structure of irreducible  module for  
$\wht{\mathfrak{gl}}_2 $. Thus,   $L_i $ possesses the structure of irreducible module  for the universal affine vertex algebra $V^1(\mathfrak{gl}_2)$ of level $1$. 
On the other hand,   the classical limit   of the  $\mathcal{V}^1(\mathfrak{gl}_2)$-module   $  \mathcal{L}_i(\mathfrak{gl}_2)$ coincides  with the $V^1(\mathfrak{gl}_2)$-module   $  L_i $. 

In the next corollary, we use the notion of irreducible module in the following sense. The topologically free $\CC[[h]]$-module V is said to be irreducible with respect to the action of associative algebra (resp. quantum vertex algebra)   if it does not possess any nontrivial topologically free $\CC[[h]]$-submodule $W$ which is invariant under the corresponding action of associative algebra (resp. quantum vertex algebra)   and satisfies that $h^n v\in W$ for $n\geqslant 1$ and $v\in V$ implies $v\in W$.

\begin{kor}\label{thm33cor2}
  The $\CC[[h]]$-modules
 $\mathcal{L}_i(\mathfrak{gl}_2)$      are irreducible  both as  modules for the quantum vertex algebra $\mathcal{V}^1(\mathfrak{gl}_2)$ and as level $1$ modules for the double Yangian $\dyg$.
\end{kor}

\begin{prf}
The irreducibility with respect to the action of $\mathcal{V}^1(\mathfrak{gl}_2)$ follows by the discussion preceding the corollary.
As for the action of the double Yangian, if the $\CC[[h]]$-submodule
$W\subset \mathcal{L}_i(\mathfrak{gl}_2)$ is  $\dyg$-invariant, by Corollary \ref{thm33cor}   we have
$$
Y_{\mathcal{L}_i}(T^-_{1n}(u_1)\ldots T^-_{n-1\ts n}(u_{n-1}),z)W\subset
(\ndo\CC^N)^{\ot (n-1)}\ot W[[z^{\pm 1},u_1,\ldots ,u_{n-1}]]
$$
for all $n\geqslant 1$. As the coefficients of matrix entries of $T^-_{1n}(u_1)\ldots T^-_{n-1\ts n}(u_{n-1})$ with $n\geqslant 1$ span an $h$-adically dense $\CC[[h]]$-submodule of $\mathcal{V}^1(\mathfrak{gl}_2)$, it is clear that $W$ is $\mathcal{V}^1(\mathfrak{gl}_2)$-invariant as well, so that   the second assertion follows.
\end{prf}

\begin{rem}
Ding--Feigin's construction of semi-infinite monomial bases    \cite{DF} 
in the case of quantum affine algebra $U_q(\wht{\mathfrak{sl}}_2)$
 suggests that the  suitable analogues  of Corollaries \ref{thm33cor} and \ref{thm33cor2} can be  established for   Etingof--Kazhdan's quantum vertex algebra \cite{EK} associated with the trigonometric $R$-matrix of type $A_1^{(1)}$. On the other hand, by using different approach, the $h$-adic quantum vertex algebra structure  over certain wide class of irreducible modules for untwisted quantum affinization algebras was recently obtained by Kong \cite{Kong}.
\end{rem}

\section*{Acknowledgement}
M.B. and S.K. would like to thank Mirko Primc for   helpful discussions.
This work has been  supported in part by Chinese National Natural Science Foundation grant nos. 12101261 and 12171303 and by Croatian Science Foundation under the project UIP-2019-04-8488.

 \linespread{1.0}

\end{document}